\documentclass{amsart}
\usepackage{amsmath,amssymb,amscd,latexsym}
\usepackage[all]{xy}

\newcommand{\Fh}{{\mathcal F}}

\newcommand{\be}{\mathbf{1}}

\newcommand{\dimnuc}{\mathrm{dim}_{\mathrm{nuc}}}

\newcommand{\id}{\mathrm{id}}

\newcounter{number}[section]

\newenvironment{nummer}{\refstepcounter{number}{\noindent\arabic{section}.\arabic{number}}}{}

\newcommand{\bn}{\noindent \begin{nummer} \rm}
\newcommand{\en}{\end{nummer}}

\newenvironment{ntheorem}{\noindent {\sc Theorem:} \it}{}
\newenvironment{nlemma}{\noindent {\sc Lemma:} \it}{}
\newenvironment{nprop}{\noindent {\sc Proposition:} \it}{}
\newenvironment{ndefn}{\noindent {\sc Definition:} \it}{}

\newenvironment{nproof}{\noindent {\sc Proof:}}{\mbox{}\hfill 
\rule[-.2ex]{.25em}{1.8ex}}

\begin{document}

\title[Minimal dynamics and $\mathrm{K}$-theoretic rigidity]{{\sc Minimal dynamics and $\mathrm{K}$-theoretic rigidity: \\ Elliott's conjecture}}

\author{Andrew S.\  Toms}
\address{Department of Mathematics and Statistics\\
York University\\ 
4700 Keele St.\\
Toronto, Ontario\\
M3J 1P3}

\email{atoms@mathstat.yorku.ca}

\author{Wilhelm Winter}
\address{School of Mathematical Sciences\\
University of Nottingham\\
University Park\\
Nottingham NG7 2RD\\ 
United Kingdom}

\email{wilhelm.winter@nottingham.ac.uk}

\date{\today}
\thanks{{\it Supported by:} EPSRC First Grant EP/G014019/1 and by an NSERC Discovery Grant}

\setcounter{section}{-1}

\begin{abstract}
Let $X$ be an infinite, compact, metrizable space of finite covering dimension and $\alpha:X \to X$ a minimal homeomorphism.  We prove that the crossed product $\mathcal{C}(X) \rtimes_\alpha \mathbb{Z}$ absorbs the Jiang-Su algebra tensorially and has finite nuclear dimension.  As a consequence, these algebras are determined up to isomorphism by their graded ordered $\mathrm{K}$-theory under the necessary condition that their projections separate traces.  This result applies, in particular, to those crossed products arising from uniquely ergodic homeomorphisms.


\end{abstract}

\maketitle

\section{Introduction}

\noindent
From its earliest days the theory of operator algebras has been entwined with dynamics, and some of the most important developments in the subject revolve around this interaction.  The group-measure space construction of Murray and von Neumann 
provided the first examples of non-type-I factors;  recently we have seen connections with orbit equivalence and associated rigidity phenomena in the remarkable works of Popa and Giordano et al. (see \cite{po} and \cite{gmps}, respectively).  

In this article we consider crossed product $\mathrm{C}^*$-algebras arising from topological dynamical systems, and prove a maximally general theorem concerning the degree to which they are determined by their graded ordered $\mathrm{K}$-theory.    


\bn
\begin{ntheorem}\label{main}
Let $\mathcal{C}$ denote the class of $\mathrm{C}^*$-algebras having the following properties:
\begin{enumerate}
\item[(i)] every $A \in \mathcal{C}$ has the form $\mathcal{C}(X) \rtimes_\alpha \mathbb{Z}$ for some infinite, compact, 
finite-dimensional, metrizable space $X$ and minimal homeomorphism $\alpha:X \to X$;
\item[(ii)] the projections of every $A \in \mathcal{C}$ separate traces.
\end{enumerate}
If $A,B \in \mathcal{C}$ and there is a graded ordered isomorphism $\phi:\mathrm{K}_*(A) \to \mathrm{K}_*(B)$, then there is
a $*$-isomorphism $\Phi:A \to B$ which induces $\phi$.
\end{ntheorem}
\en

This result was conjectured by G.\ A.\ Elliott in 1990 as part of his wider program to classify separable amenable $\mathrm{C}^*$-algebras.  The hypotheses of minimality for $\alpha$, finite-dimensionality for $X$, and the separation of traces by projections are all known to be necessary;  the necessity of finite-dimensionality for $X$ was established recently by Giol and Kerr in \cite{gk}.   If one imposes unique ergodicity on $\alpha$, then condition (ii) is unnecessary.   
Our result is the culmination of a sequence of earlier important results due to Elliott and Evans (\cite{ee}), H. Lin and Phillips (\cite{lp}), and the second named author (\cite{w2, w3}).  It covers irrational rotation algebras (which have many projections) as well as the (projectionless) $\mathrm{C}^{*}$-algebras associated to minimal homeomorphisms of odd spheres considered in \cite{Con:Thom}.

Our proof does not employ a decomposition for $\mathcal{C}(X) \rtimes_\alpha \mathbb{Z}$ as an inductive limit of type I $\mathrm{C}^*$-algebras;  rather, it implies that one exists.  Using \cite{StrWin:TAS}, we reduce the problem to a setting similar to the one considered in \cite{lp}, where an abundance of projections was assumed.
This reduction is made by applying a general classification result developed by the second named author in \cite{w3}, the key hypothesis of which is the condition that the $\mathrm{C}^*$-algebras considered all absorb the Jiang-Su algebra $\mathcal{Z}$ tensorially.  
(This condition is necessary for $\mathrm{K}$-theoretic rigidity results akin to Theorem \ref{main}, both in general and in the case of tracial algebras (see \cite{r2} and \cite{t}, respectively).  We refer the reader to \cite{et} and to \cite{RorWin:Z-revisited} for a complete discussion of the Jiang-Su algebra  and its relevance to Elliott's program.)  The bulk of our effort is concentrated on proving the following result.

\bn
\begin{ntheorem}\label{zstab}
Let $X$ be an infinite, compact, finite-dimensional metrizable space and \mbox{$\alpha:X \to X$} a minimal homeomorphism.  It follows that
\[
(\mathcal{C}(X) \rtimes_\alpha \mathbb{Z}) \otimes \mathcal{Z} \cong \mathcal{C}(X) \rtimes_\alpha \mathbb{Z}.
\]
\end{ntheorem}
Notice that we do not require projections to separate traces in this theorem.  It is conjectured that Theorem \ref{main} continues to hold in the absence of condition (ii), provided that one augments the invariant $\mathrm{K}_*$ by the simplex of tracial states (here identified with the $\alpha$-invariant Borel probability measures on $X$).  We expect that Theorem \ref{zstab} will prove crucial to the solution of this conjecture, too.  
\en

The techniques that we develop in proving Theorem \ref{zstab} also allow us to bound the nuclear dimension of the crossed products we consider.  This dimension, based on completely positive approximations of the identity map, generalizes the usual covering dimension of a locally compact Hausdorff space to the realm of nuclear $\mathrm{C}^*$-algebras.  Its close cousin, the decomposition rank, has already proved to be a very powerful tool in efforts to further Elliott's classification program, and there is much evidence to suggest that the nuclear dimension will be similarly important.  The key question is this:  when is the nuclear dimension finite?

\bn
\begin{ntheorem}\label{dimnuc}
Let $X$ be an infinite, compact, finite-dimensional metrizable space and \mbox{$\alpha:X \to X$} a minimal homeomorphism.  It follows that
the nuclear dimension of $\mathcal{C}(X) \rtimes_\alpha \mathbb{Z}$ is at most $2 \mathrm{dim}(X)+1$.
\end{ntheorem}

We note that the $\mathrm{C}^*$-algebras of Theorem \ref{main} have nuclear dimension at most $2$, and that the same is most likely true of the algebras considered in Theorem \ref{dimnuc}.  This improved bound, however, relies on the existence of a special inductive limit decomposition for the crossed product, and this, in turn, relies on the classification theorem itself.
\en

Our paper is organised as follows:  Section 1 collects some preliminary results, in Section 2 we establish simplicity and finite decomposition rank for some natural subalgebras of $\mathcal{C}(X) \rtimes_\alpha \mathbb{Z}$, and in Sections 3, 4, and 5 we prove Theorems \ref{dimnuc}, \ref{zstab}, and \ref{main}, respectively.

\newpage

\section{Preliminaries}

\bn
\label{function-commutator}
\begin{nprop}
Let $\mathcal{G} \subset \mathcal{C}([0,1])$ be a finite subset of positive functions. Then, for any $\eta >0$ there is $\delta >0$ such that the following holds:

If $A$ is a $\mathrm{C}^{*}$-algebra and $h \in A$ is a positive element of norm at most one, then if $b \in A$ satisfies 
\[
\|[b,h]\| < \delta,
\]
this implies 
\[
\|[b,f(h)]\| < \eta
\]
for all $f \in \mathcal{G}$. 
\end{nprop}

\begin{nproof}
It will suffice, by uniform density, to consider the case where $\mathcal{G}$ is a finite set of polynomials.  Let $k$ be the largest degree of any
$f \in \mathcal{G}$.  Observe that if $\|[b,h]\| < \delta$, then for any natural number $i \leq k$ we have
\[
bh^i \approx_\delta hbh^{i-1} \approx_\delta h^2bh^{i-2} \approx_\delta \cdots \approx_\delta h^i b,
\]
where $\approx_\delta$ denotes the relation of being at norm distance strictly less than $\delta$.
It follows that $\|[b,h^i]\| < i \delta \leq k \delta$.  

Let $f(x) = a_0+a_1x+ \cdots a_k x^k$ be a generic element of $\mathcal{G}$, and let $\eta>0$ be given.  We set $M_f = \max\{1,|a_0|,|a_1|,\ldots,|a_k| \}$ and 
\[
M = \max_{ f(x) \in \mathcal{G} } M_f.
\]
With $\delta = \eta/(k^2M)$ we compute:
\begin{eqnarray*}
\|[b,f(h)] \|& = & \left\| \sum_{i=0}^k   [b,a_i h^i] \right \| \\
& \leq & M \left \| \sum_{i=0}^k [b,h^i] \right \| \\
& \leq & M \sum_{i=0}^k \| [b,h^i] \| \\
& < & Mk(k \delta) \\
& \leq & \eta
\end{eqnarray*}

\end{nproof}
\en

\bn
\label{C-gamma}
\begin{nprop}
Let $C$ be a $\mathrm{C}^{*}$-algebra, $\mathcal{G} \subset C$ a self-adjoint subset generating $C$ as a $\mathrm{C}^{*}$-algebra, $\mathcal{E} \subset C$ finite and $\theta > 0$. Then, there are $\beta >0$ and a finite set $\mathcal{G}' \subset \mathcal{G}$ such that the following holds:

If $B$ is another $\mathrm{C}^{*}$-algebra and $\gamma: C \to B$ is a $*$-homomorphism, then if $b \in B$ has norm at most one and satisfies 
\[
\|[b,\gamma(g)]\| < \beta
\]
for $g \in \mathcal{G}'$, we have
\[
\|[b,\gamma(e)]\| < \theta
\]
for $e \in \mathcal{E}$.
\end{nprop}

\begin{nproof}
Let $\theta>0$ be given.  Since $\mathcal{G}$ generates $C$, there is a natural number $k$ such that each $e \in \mathcal{E}$ is at norm distance strictly less than $\theta/4$ from a sum of at most $k$ monomials of the form
\[
\alpha x_1 x_2 \cdots x_l,
\]
where $x_j \in \mathcal{G}$, $\alpha \in \mathbb{C}$, and $l \leq k$.  We denote this sum by $g_e$, so
that $\| e-g_e\| < \theta/4$, and we let $\mathcal{G}'$ denote the set of all elements of $\mathcal{G}$ which occur in a monomial of $g_e$ as $e$ ranges over $\mathcal{E}$.  Now, for any $\mathrm{C}^*$-algebra $B$, $b \in B$ of norm at most one and $*$-homomorphism $\gamma:C \to B$ we have
\begin{eqnarray*}
\|[b,\gamma(e)]\| & = & \|[b,\gamma(g_e) - \gamma(e-g_e)]\| \\
& \leq & \|[b,\gamma(g_e)]\| + \|[b,\gamma(e-g_e)]\| \\
& \leq & \|[b,\gamma(g_e)]\| + 2\|b\| \cdot \|\gamma(e-g_e)\| \\
& < & \|[b,\gamma(g_e)]\| + \theta/2.
\end{eqnarray*}
Thus, to complete the proof of the Proposition, we must show that there is $\beta>0$ such that for any $B$, $b$ and $\gamma$ as above we have the following statement:  if
\[
\|[b,\gamma(g)]\| < \beta
\]
for every $g \in \mathcal{G}'$, then
\[
\|[b,\gamma(g_e)]\| < \theta/2
\]
for every $e \in \mathcal{E}$.  Since $\|[b,g_e/\eta]\| < \theta/(2\eta)$ implies $\|[b,g_e]\| < \theta/2$ for any $\eta >0$, we may scale the $g_e$s and assume that the elements of $\mathcal{G}'$ have norm at most one.

Let $\gamma:C \to B$ be given, and let 
\[
\alpha x_1 x_2 \cdots x_l
\]
be a monomial appearing in the sum which constitutes some $g_e$.  Set $y_j = \gamma(x_j)$, and note that $\|y_j\| \leq 1$.  The assumption $\|[b,y]\|< \beta$ for
every $y \in \gamma(\mathcal{G}')$ then yields
\begin{eqnarray*}
\lefteqn{\| (y_1 y_2 \cdots y_{j-1} b y_j \cdots y_l) - (y_1 y_2 \cdots y_j b y_{j+1} \cdots y_l) \| } \\
& = & \| (y_1 \cdots y_j)[b,y_j](y_{j+1} \cdots y_l) \|  \\
& \leq & \| [b,y_j] \| \\
& < & \beta.
\end{eqnarray*}
It follows that
\begin{eqnarray*}
y_1 y_2 \cdots y_l b & \approx_\beta  &  y_1 y_2 \cdots y_{l-1} b y_l \\
& \approx_\beta & y_1 y_2 \cdots y_{l-2} b y_{l-1} y_l \\
& \cdots & \\
& \approx_\beta & b  y_1 y_2 \cdots y_l,
\end{eqnarray*}
and so
\[
\| [b,  y_1 y_2 \cdots y_l] \| < k \beta.
\]

We have
\[
\gamma(g_e) = \sum_{j=1}^k \alpha_j y_{j,1} y_{j,2} \cdots y_{j,l_j},
\]
where $\alpha_j \in \mathbb{C}$ and $y_{j,t} \in \gamma(\mathcal{G}')$.  (Recall that $k$ does not depend on $e$, and that $l_j \leq k$ regardless of $i$.)  Set $M_e = \max \{1,|\alpha_1|,\ldots,| \alpha_{k} | \}$, and 
\[
M = \max_{e \in \mathcal{E}} M_e
\]
Set 
\[
\beta = \theta/(2k^2M),
\]
and observe that $\beta$ does not depend on $\gamma$.  Now, we compute
\begin{eqnarray*}
\| [b,\gamma(g_e)] \| & = & \left \| \sum_{j=1}^k [b, \alpha_j  y_{j,1} y_{j,2} \cdots y_{j,l_j}] \right \| \\
& \leq & M \left \| \sum_{j=1}^k [b,y_{j,1} y_{j,2} \cdots y_{j,l_j}] \right \| \\
& \leq & M \sum_{j=1}^k \| [b,y_{j,1} y_{j,2} \cdots y_{j,l_j}] \| \\
& < & Mk(k \beta)/2 = \theta/2,
\end{eqnarray*}
as required.
\end{nproof}
\en

\section{Local simplicity}

\bn
\label{A-Y-def}
The arguments in this and the next section will employ the notion of a recursive subhomogeneous $\mathrm{C}^*$-algebra. We refer the reader to \cite{p2} for the basic definitions and terminology related to these algebras.  We will also use heavily certain natural subalgebras of crossed product $\mathrm{C}^*$-algebras, which arise as follows.  Let $X$ be a compact metrizable space, $\alpha: X \to X$ a homeomorphism, $Y$ a closed subset of $X$, and $u$ the unitary implementing the action of $\alpha$ in $A = \mathcal{C}(X) \rtimes_\alpha \mathbb{Z}$;  we define
\begin{equation}\label{ay}
A_Y = \mathrm{C}^*(\mathcal{C}(X), u \mathcal{C}_0(X \backslash Y)).
\end{equation}
\en

\bn
\label{p-local-simplicity}
\begin{nprop}
Let $X$ be an infinite, compact, metrizable, finite-dimensional space with a minimal homeomorphism $\alpha$. If $x_{0}, \, x_{1} \in X$ have
disjoint orbits under $\alpha$, then the subalgebras $A_{\{x_{0}\}}$,  $A_{\{x_{0},x_{1}\}}$ and $A_{\{x_{1}\}}$ are all simple and have finite 
decomposition rank.
\end{nprop}

\begin{nproof}
Let us first prove the statement about finite decomposition rank.  Use $Y$ to denote any of $\{x_0\}$, $\{x_1\}$, or $\{x_0,x_1\}$.  Let $Y_m$ 
be a decreasing sequence of closed subsets of $X$ such that $\mathrm{int}(Y_m) \neq \emptyset$ and $\cap_{m=1}^\infty Y_m = Y$.  It is immediate
that 
\[
A_Y = \mathrm{C}^*(\mathcal{C}(X),u \mathcal{C}_0(X \backslash Y) ) = \lim_{m \to \infty} \mathrm{C}^*(\mathcal{C}_0(X),u \mathcal{C}(X \backslash Y_m)) 
= \lim_{m \to \infty} A_{Y_m}
\]
(see (\ref{ay})).  It is proved in \cite{qlp} 
that $A_{Y_m}$ is a recursive subhomogeneous $\mathrm{C}^*$-algebra of topological dimension
at most $\mathrm{dim}(X)$.  It then follows from the main theorem of \cite{w4} that each $A_{Y_m}$ has decomposition rank at most 
$\mathrm{dim}(X)$.  Decomposition rank is lower semicontinuous with respect to inductive limits, whence $A_Y$ has decomposition rank at most
$\mathrm{dim}(X)$, too (see \cite{kw}).

The simplicity of $A_Y$ is established by \cite[Proposition 2.5]{lp} in the case that $Y$ is a singleton.  We adapt their proof to 
establish the simplicity of $A_{\{x_0,x_1\}}$.  Let $I \subset A_{\{x_0,x_1\}}$ be a nonzero ideal.  It follows that $\mathcal{C}(X) \cap I$ is an
ideal of $\mathcal{C}(X)$, and therefore has the form $\mathcal{C}_0(U)$ for some open subset $U$.  Explicitly,
\[
U = \{ x \in X \ | \ \exists f \in \mathcal{C}(X) \cap I : f(x) \neq 0 \}.
\]
The proof of \cite[Proposition 2.5]{lp} shows that $U \neq \emptyset$.  

We claim that $\alpha^{-1} ( U \backslash \{\alpha(x_0),\alpha(x_1)\}) \subset U$.  To see this, fix $z \in U \backslash \{\alpha(x_0),\alpha(x_1)\}$. 
Choose $f \in \mathcal{C}(X) \cap I$ such that $f(z) \neq 0$, and choose $g \in \mathcal{C}_0(X \backslash \{x_0,x_1\})$ such that
$g(\alpha^{-1}(z)) \neq 0$.  Now $ug \in A_{\{x_0,x_1\}}$, and so
\[
(ug)^* f (ug) = \overline{g} u^* f u g = |g|^2(f \circ \alpha)
\]
is in $\mathcal{C}(X) \cap I$.  Since $|g|^2(f \circ \alpha)$ is not zero at $\alpha^{-1}(z)$, we see that $\alpha^{-1}(z) \in U$.

We also claim that $\alpha(U \backslash \{x_0,x_1\}) \subset U$.  Let $z \in U \backslash \{x_0,x_1\}$, and choose $f \in \mathcal{C}(X) \cap I$ and 
$g \in \mathcal{C}_0(X \backslash \{x_0,x_1\})$ to be nonzero at $z$.  Now $ug \in A_{\{x_0,x_1\}}$, and so
\[
(ug) f (ug)^* = u(|g|^2f)u^* = (|g|^2f) \circ \alpha^{-1}
\]
is in $\mathcal{C}(X) \cap I$.  Since $(|g|^2f) \circ \alpha^{-1}$ is not zero at $\alpha(z)$, we see that $\alpha(z) \in U$.

Set $Z = X \backslash U$.  We will prove that $Z = \emptyset$ by contradiction.

First suppose that $z \in Z$ and $z \notin \mathrm{Orb}(x_0) \cup \mathrm{Orb}(x_1)$.    If $k>0$ then, since $\alpha^k(z) \neq \alpha(x_0), \alpha(x_1)$, we use the first of the two
claims above to conclude that $\alpha^{k-1}(z) \in U \backslash \{\alpha(x_0),\alpha(x_1)\}$.  Iterating this procedure $k$ times yields $z \in U$, 
contradicting $z \in Z$.  We conclude that $\alpha^k(z) \in Z$ for any $k>0$.  
Similarly,  for  $k < 0$,  since $\alpha^k(z) \neq x_0,x_1$ we can use the second of the two claims above to 
conclude that $\alpha^{k+1}(z) \in U \backslash \{x_0,x_1\}$.  Iterating this procedure $k$ times yields $z \in U$, contradicting $z \in Z$.
We conclude again that $\alpha^k(z) \in Z$, this time for any $k<0$.
We have proved that $\mathrm{Orb}(z) \subset Z$, but this contradicts the minimality of $\alpha$.  We conclude that if  
$z \in Z$, then $z \in \mathrm{Orb}(x_0) \cup \mathrm{Orb}(x_1)$.

Suppose, without loss of generality, that $z = \alpha^k(x_0)$.  Assume that $k>0$, and that $\alpha^n(x_0) \in U$
for some $n > k$.   Using $\mathrm{Orb}(x_0) \cap \mathrm{Orb}(x_1) = \emptyset$ and $n>1$ we see that in fact
$\alpha^n(x_0) \in U \backslash \{\alpha(x_0),\alpha(x_1) \}$.  Applying the first of our two claims above, we see that 
$\alpha^{n-1}(x_0) \in U \backslash \{\alpha(x_0),\alpha_(x_1)\}$.  Iterating this argument yields $\alpha^k(x_0) \in U$, a contradiction, so we must have
$\{\alpha^n(x_0) \ | \ n \geq k \} \subset Z$.  This, however, violates minimality.  Similarly, suppose that $k \leq 0$, and that 
$\alpha^n(x_0) \in U $ for some $n < k$.  Using $\mathrm{Orb}(x_0) \cap \mathrm{Orb}(x_1) = \emptyset$ and $n < 0$ we see that
in fact $\alpha^n(x_0) \in U \backslash \{x_0,x_1\}$.  Applying the 
second of our two claims above, we see that $\alpha^{n+1}(x_0) \in U \backslash \{ x_0,x_1\}$.  Iterating this argument yields $\alpha^k(x_0) \in U$,
a contradiction, so we conclude that $\{ \alpha^n(x_0) \ | \ n \leq k\} \subset Z$.  This again violates minimality.  

Since every possible choice of $z \in Z$ leads to a contradiction, we conclude that $Z = \emptyset$, $U = X$, and so $1_A \in I$.  It follows that
$A_{\{x_0,x_1\}}$ is simple, as desired.
\end{nproof}
\en

\section{Nuclear dimension}

\noindent
In this section we recall the notion of nuclear dimension as introduced in \cite{wz} and prove our Theorem~\ref{dimnuc}.

\bn
\label{dimnuc-def}
\begin{ndefn}
A  $\mathrm{C}^*$-algebra $A$ has nuclear dimension at most $n$, $\dimnuc A \le n$, if there exists a net 
$(F_{\lambda},\psi_{\lambda},\varphi_{\lambda})$ of finite-dimensional c.p.\ approximations for $A$ (i.e., $F_{\lambda}$ are finite dimensional $\mathrm{C}^{*}$-algebras, and $\psi_{\lambda}: A \to F_{\lambda}$  and  
$\varphi_{\lambda}: F_{\lambda} \to A$ are completely positive maps  for all $\lambda$) such that 
\begin{enumerate}
\item[(i)] $\varphi_{\lambda} \circ \psi_{\lambda} (a) \to a$ for any  $a \in A$
\item[(ii)] $\|\psi_{\lambda}\| \le 1$ for each $\lambda$
\item[(iii)] each $F_{\lambda}$ decomposes into $n+1$ ideals, $F_{\lambda}= F_{\lambda}^{(0)} \oplus \ldots \oplus F_{\lambda}^{(n)}$,  
such that $\varphi_{\lambda}|_{F_{\lambda}^{(i)}} $ is an order zero contraction (i.e., it is c.p.c.\ and preserves orthogonality) for $i=0, \ldots ,n$.
\end{enumerate}
\end{ndefn}
\en

\bn
\label{almost-central-h}
\begin{nprop}
Let $X$ be an infinite, compact, metrizable space with a minimal homeomorphism $\alpha$. It follows that for any $\delta >0$ and $\mathcal{F} \subset \mathcal{C}(X) \rtimes_{\alpha} \mathbb{Z}$ finite, there are a positive normalized element $h \in \mathcal{C}(X)$ and distinct points $x_0,x_1 \in X$ such that
\[
h(x_{0})=0, \, h(x_{1})=1\]
and 
\[
\|[h,b]\| < \delta
\]
for all $b \in \mathcal{F}$.  We may moreover arrange that $\mathrm{Orb(x_0)} \cap \mathrm{Orb}(x_1) = \emptyset$.
\end{nprop}

\begin{nproof}
Let $u$ be the unitary operator in $A := \mathcal{C}(X) \rtimes_\alpha \mathbb{Z}$ implementing $\alpha$.
Applying Proposition \ref{C-gamma} with $C=A$, $\gamma = \mathrm{id}_A$, $\mathcal{G} = \{\mathcal{C}(X),u,u^*\}$, $\mathcal{E} = \Fh$, and $\theta = \delta$, 
we see that there is $0 < \eta < \delta$ 
such that if $h \in \mathcal{C}(X)$ is positive with the property that 
\[
\|[h,u^*]\|=\|[h,u]\|<\eta,
\]
then $\|[h,b]\|< \delta$ for every $b \in \Fh$.  Thus, to establish the commutator estimate of this Proposition, 
we need only show that $h$ can be chosen to commute with $u$ to within $\eta$.

Set $n = \lceil 1/\eta \rceil + 1$, and fix some $x_0 \in X$.  
By minimality and the fact that $X$ is infinite, we can find an open neighbourhood $U$ of $x_0$ such that the 
sets $\alpha^i(U)$ are pairwise disjoint for $i \in \{0,1,\ldots,2n\}$.  Let $f \in \mathcal{C}_0(U)$ be a positive function 
of norm one with the property that $f \equiv 1$ on an open neighbourhood $V$ of $x_0$ such that $\overline{V} \subset U$.  
Now define $h \in \mathcal{C}_0(U \cup \alpha(U) \cup \cdots \cup \alpha^{2n}(U))$ via the formula
\[
h|_{\alpha^i(U)} = \left\{ \begin{array}{ll} (i/n)(f \circ \alpha^{-i}), & i \leq n \\ (2-i/n) (f \circ \alpha^{-i}), & n < i \leq 2n \end{array} \right. .
\]
Our assumptions on $X$ and $\alpha$ imply that there is more than one orbit under $\alpha$.  By minimality we can
find $x_1 \in \alpha^n(V)$ such that $\mathrm{Orb(x_0)} \cap \mathrm{Orb(x_1)} = \emptyset$, whence $h(x_1)=1$ and $h(x_0)=0$.  
Clearly, $h$ is positive and of norm one.

It remains to prove that $\| uh - hu \|< \eta$.  This is equivalent to showing that 
\[
\| uhu^* - h \| < \eta.
\]
Now 
\[
uhu^* = h \circ \alpha^{-1} \in \mathcal{C}_0(U \cup \alpha(U) \cup \cdots \cup \alpha^{2n}(U))
\]
is given by the formula
\[
uhu^*|_{\alpha^i(U)} = \left\{ \begin{array}{ll} 0, & i = 2n \\ (2-(i+1)/n)(f \circ \alpha^{-i}), & n-1 < i \leq 2n-1 \\ ((i+1)/n)(f \circ \alpha^{-i}), & 0 \leq i \leq n-1 \end{array}
\right. ,
\] 
and comparing this with the formula for $h$ we see that
\[
\|uhu^*-h\| = (1/n)\|f\| = 1/n < \eta,
\]
as required.
\end{nproof}
\en

\bn
\begin{ntheorem}
Let $X$ be an infinite, compact, metrizable, finite-dimensional space with a minimal homeomorphism $\alpha$. Then,
\[
\dimnuc (\mathcal{C}(X) \rtimes_{\alpha} \mathbb{Z}) \le 2 \dim X + 1.
\]
\end{ntheorem}

\begin{nproof}
Set $A = \mathcal{C}(X) \rtimes_{\alpha} \mathbb{Z}$, and let there be given a finite subset $\Fh = \{a_1,\ldots,a_n\}$ of $A$ and a 
tolerance $\epsilon>0$.  We must find a c.p. approximation for $\Fh$ to within $\epsilon$ which has the form described in Definition
\ref{dimnuc-def}. 

Since $\mathcal{C}(X)$ and the implementing 
unitary $u$ together generate $A$, there is a finite set $\Fh' = \{a_1',\ldots,a_n'\}$ with the following properties:
\begin{enumerate}
\item[(i)] $\|a_i-a_i'\| < \epsilon/24$, for each $i \in \{1,\ldots,n\}$;
\item[(ii)] there exists $k \in \mathbb{N}$ such that each $a_i'$ is the sum of at most $k$ monomials of the form
\[
f_1 u^{j_1} f_2 u^{j_2} \cdots f_k u^{j_k}, \ j_l \in \{-1,0,1\}, \ f_l \in \mathcal{C}(X), \ 1 \leq l \leq k.
\]
\end{enumerate}

Apply Proposition \ref{almost-central-h} to the singleton $\{u\}$ and a value of $\delta$ so small that the resulting function $h$ has the property that
\[
\|[u^*,h^{i/k}]\|=\|[u,h^{i/k}]\| < \epsilon/(24k^2)
\]
and
\[
\|[u^*,(1-h)^{i/k}]\|=\|[u,(1-h)^{i/k}]\| < \epsilon/(24k^2)
\]
for each $i \in \{1,\ldots,k\}$  (this is possible by Proposition \ref{function-commutator}).  Now for any monomial as in (ii) above, we have
\begin{equation}\label{monoest}
\|(f_1 u^{j_1} f_2 u^{j_2} \cdots f_k u^{j_k})h -  (f_1 u^{j_1} h^{1/k} f_2 h^{1/k} u^{j_2} \cdots f_k u^{j_k}h^{1/k})\| < \epsilon/24k, 
\end{equation}
where we have inserted $h^{1/k}$ to the left of $u^{j_l}$ if $j_l = -1$, and to the right if $j_l = 0,1$.  To see why this is so, 
observe that we may pass from
\begin{equation}\label{mono}
(f_1 u^{j_1} f_2 u^{j_2} \cdots f_k u^{j_k})h 
\end{equation}
to 
\begin{equation}\label{mono2}
(f_1 u^{j_1} h^{1/k} f_2 h^{1/k} u^{j_2} \cdots f_k u^{j_k}h^{1/k})
\end{equation}
in at most $k$ steps, each of which involves passing a function of the form $h^{l/k}$ from one side of a $u$ or $u^*$ to the other.
Keep in mind that $h^{l/k}$ commutes with each $f_{t}$ for all sensible $l$ and $t$.
Each such step has a cost in norm which is less than or equal to  
\[
\|[u^*,h^{i/k}]\|=\|[u,h^{i/k}]\| < \epsilon/(24k^2).
\]

Each of $f_l$, $h^{1/k}u^{-1}$,
and $u h^{1/k}$ is contained in 
\[
A_{\{x_0\}} = \mathrm{C}^*(\mathcal{C}(X), u \mathcal{C}_0(X \backslash \{x_0\}).
\]
The monomial of (\ref{mono2}) is a product of such elements, hence also in $A_{\{x_0\}}$.  Since each $a_i'h$ is a sum of at most
$k$ monomials as in (\ref{mono}) above, and since each such monomial is at distance strictly less that $\epsilon/(24k)$ from a monomial as 
in (\ref{mono2}), we conclude that $a_i'h$ is at distance at most $\epsilon/24$ from an element $b_i^{(0)}$ of $A_{\{x_0\}}$.  
Set $\mathcal{G}_0=\{b_1^{(0)},\ldots,b_n^{(0)}\}$.  An argument similar to the one just presented shows that $a_i'(1-h)$ is at distance
at most $\epsilon/24$ from some
\[
b_i^{(1)} \in A_{\{x_1\}} = \mathrm{C}^*(\mathcal{C}(X),u \mathcal{C}_0(X \backslash \{x_1\})).
\]
Set $\mathcal{G}_1 = \{b_1^{(1)},\ldots,b_n^{(1)}\}$.  Using the fact that $\|a_i - a_i'\| < \epsilon/24$, and also that we may assume that $h$ commutes
with each $a_i$ to within an arbitrarily small tolerance, we conclude that
\[
\| h^{1/2} a_i h^{1/2} - b_i^{(0)} \|, \ \|(1-h)^{1/2} a_i (1-h)^{1/2} - b_i^{(1)} \| < \epsilon/12.
\]

From Proposition \ref{p-local-simplicity}, we know that both $A_{\{x_0\}}$ and $A_{\{x_1\}}$ have decomposition rank $d = \mathrm{dim}(X) < \infty$.
We can therefore find finite-dimensional $\mathrm{C}^*$-algebras $F^{(0)}$ and $F^{(1)}$ and c.p. contractions $\phi^{(j)}:A_{\{x_j\}} \to F^{(j)}$ and
$\psi^{(j)}:F^{(j)} \to A_{\{x_j\}}$ with the property that
\[
\|\psi^{(j)} \circ \phi^{(j)}(b^{(j)}_i) - b_i^{(j)}\| < \epsilon/12, \ j \in \{0,1\}, \ i \in \{1,\ldots,n\}.
\]
We may moreover take the $\psi^{(j)}$ to be $d$-decomposable in the sense of \cite{kw}.  By Arveson's Extension Theorem, the $\phi^{(j)}$ may be extended to 
c.p. contractions 
\[
\overline{\phi^{(j)}}: A \to F^{(j)}, \ j = 0,1, 
\]
and we can thus define a c.p. contraction $\overline{\phi}:A \to F^{(0)} \oplus F^{(1)}$ by the formula
\[
\overline{\phi}(x) = \overline{\phi^{(0)}}\left(h^{1/2}xh^{1/2}\right) \oplus \overline{\phi^{(1)}}\left((1-h)^{1/2}x(1-h)^{1/2}\right).
\]
We also define a c.p. map $\overline{\psi}:F^{(0)} \oplus F^{(1)} \to A$ by the formula
\[
\overline{\psi}(x) = \psi^{(0)}(x) + \psi^{(1)}(x),
\]
assuming that $\psi^{(0)}|_{F^{(1)}} \equiv 0$ and $\psi^{(1)}|_{F^{(0)}} \equiv 0$.
By the $d$-decomposability of the $\psi^{(j)}$, we can write
\[
F^{(j)} = F^{(j)}_0 \oplus F^{(j)}_1 \oplus \cdots \oplus F^{(j)}_d
\]
so that $\psi^{(j)}|_{F^{(j)}_l}$ is an order zero contraction for each $l \in \{0,1,\ldots,d\}$ and $j \in \{0,1\}$.  
It follows immediately that the restriction of $\overline{\psi}$ to any of the direct summands of
\[
F^{(0)}_0 \oplus F^{(0)}_1 \oplus \cdots \oplus F^{(0)}_d \oplus F^{(1)}_0 \oplus F^{(1)}_1 \oplus \cdots \oplus F^{(1)}_d (= F^{(0)} \oplus F^{(1)})
\]
is an order zero contraction, so that $\overline{\psi}$ is $(2d+1)$-decomposable. 

Now we estimate
\begin{eqnarray*}
\lefteqn{ \|\overline{\psi} \circ \overline{\phi} (a_i) - a_i \| }\\
& = & \| \psi^{(0)} \circ \overline{\phi^{(0)}}(h^{1/2} a_i h^{1/2}) + \psi^{(1)} \circ \overline{\phi^{(1)}}((1-h)^{1/2}a_i(1-h)^{1/2}) - a_i \| \\
& < & \| \psi^{(0)} \circ \phi^{(0)}(b_i^{(0)}) + \psi^{(1)} \circ \phi^{(1)}(b_i^{(1)}) - a_i \| + \epsilon/6 \\
& \leq & \| b_i^{(0)} + b_i^{(1)} - a_i \| + \epsilon/3 \\
& < & \| a_i'h + a_i'(1-h) - a_i \| + \epsilon/2 \\
& = & \| a_i' - a_i \| + \epsilon/2 \\
& < & \epsilon.
\end{eqnarray*}
Since $\|\overline{\phi}\| \leq 1$ and $\overline{\psi}$ is $(2d+1)$-decomposable, we conclude that the nuclear dimension of $A$ is at most
$2d+1$, as desired.
\end{nproof}
\en

\section{$\mathcal{Z}$-stability}

In this section we combine ideas from \cite{DadWin:trivial-fields} and \cite{w2} with those of the preceding section to prove our Theorem~\ref{zstab}.

\bn
\label{almost-central-D-stable}
The following is only a minor modification of \cite[Theorem~7.2.2]{Ror:encyc}. See \cite{tw:ssa} for an introduction to strongly self-absorbing $\mathrm{C}^{*}$-algebras.

\begin{nprop}
Let $A$ and $\mathcal{D}$ be separable $\mathrm{C}^{*}$-algebras, $\mathcal{D}$ strongly self-absorbing. Suppose that for any $\epsilon >0$ and finite subsets $\mathcal{F} \subset A$ and $\mathcal{E} \subset \mathcal{D}$ there is a unital $*$-homomorphism 
\[
\zeta: \mathcal{D} \to \mathcal{M}(A)_{\infty} := \prod_{\mathbb{N}} \mathcal{M}(A) / \bigoplus_{\mathbb{N}} \mathcal{M}(A)
\]
such that
\[
\| [b,\zeta(z)] \| < \epsilon
\]
for $b \in \mathcal{F}$ and $z \in \mathcal{E}$. Then, $A$ is $\mathcal{D}$-stable.
\end{nprop}

\begin{nproof}
From the hypotheses we obtain a sequence of unital $*$-homomorphisms  $\zeta_{n}: \mathcal{D} \to \mathcal{M}(A)_{\infty}$ such that 
\[
[\zeta_{n}(z),a] \stackrel{n \to \infty}{\longrightarrow} 0 \; \forall \,  z \in \mathcal{D}, \, a \in A.
\]
As $\mathcal{D}$ is nuclear, we can lift each $\zeta_{n}$ to  a u.p.c.\ map
\[
\bar{\zeta}_{n}: \mathcal{D} \to \prod_{\mathbb{N}} \mathcal{M}(A).
\]
These give rise to a u.p.c.\ map
\[
\bar{\zeta}:\mathcal{D} \to \prod_{\mathbb{N} \times \mathbb{N}} \mathcal{M}(A).
\]
Using separability of $\mathcal{D}$ and $A$, it is now straightforward to construct a diagonal sequence of approximately multiplicative u.c.p.\ maps
\[
\bar{\zeta}_{n_{k},m_{k}}: \mathcal{D} \to \mathcal{M}(A)
\]
such that
\[
[\bar{\zeta}_{n_{k},m_{k}}(z), a] \stackrel{k \to \infty}{\longrightarrow} 0 \; \forall \, z \in \mathcal{D},\, a \in A.
\]
These induce a unital $*$-homomorphism 
\[
\tilde{\zeta}: \mathcal{D} \to \mathcal{M}(A)_{\infty} \cap A',
\]
which in turn means that $A$ is $\mathcal{D}$-stable, cf.\ \cite[Theorem~7.2.2]{Ror:encyc}.
\end{nproof}

We will apply  Proposition \ref{almost-central-D-stable} in the case $\mathcal{D} = \mathcal{Z}$, where $\mathcal{Z}$
is the Jiang-Su algebra, cf.\ \cite{tw:ssa}.
\en

\bn
\label{diagonally-central}
\begin{nprop}
Let $A$ be a unital $\mathcal{Z}$-stable $\mathrm{C}^{*}$-algebra. If 
\[
B \subset A_{\infty}
\]
is a separable subset, there is a $*$-homomorphism 
\[
\varrho: A \otimes \mathcal{Z} \to A_{\infty}
\]
such that 
\[
\varrho|_{A \otimes \be_{\mathcal{Z}}}  = \iota_{A}
\]
(where $\iota_{A}$ denotes the canonical embedding of $A$ into $A_{\infty}$) and
\[
[\varrho(\be_{A} \otimes \mathcal{Z}), B]=0.
\]
\end{nprop}

\begin{nproof}
Since $A$ is unital and $\mathcal{Z}$-stable, by \cite[Theorem~7.2.2]{Ror:encyc} there is a unital $*$-homomorphism 
\[
\tilde{\varrho}: \mathcal{Z} \to A_{\infty} \cap A'.
\]
Now by \cite[Lemma~4.5]{HirWin:rokhlin}, there is a unital $*$-homomorphism
\[
\bar{\varrho}: \mathcal{Z} \to A_{\infty} \cap A' \cap C'
\]
(in \cite[Lemma~4.5]{HirWin:rokhlin}, the subspace $B$ was assumed to be a $\mathrm{C}^{*}$-subalgebra of $A_{\infty} \cap A'$, but inspection of the proof shows that the result remains valid if $B$ is just a subset of $A_{\infty}$).
\end{nproof}
\en

\bn
\label{C-Z-trivial}
\begin{nprop}
Let $\{d_{0},d_{1/2},d_{1}\}$ be a partition of unity for $[0,1]$. Then,
\[
C:= \mathrm{C}^{*}(d_{0} \otimes \mathcal{Z} \otimes \be_{\mathcal{Z}} \otimes \be_{\mathcal{Z}} \cup d_{1/2} \otimes \be_{\mathcal{Z}} \otimes \mathcal{Z} \otimes \be_{\mathcal{Z}} \cup d_{1} \otimes \be_{\mathcal{Z}} \otimes \be_{\mathcal{Z}} \otimes \mathcal{Z} ) \subset \mathcal{C}([0,1]) \otimes \mathcal{Z} \otimes \mathcal{Z} \otimes \mathcal{Z}
\]
is isomorphic to $\mathrm{C}^{*}(d_{0},d_{1/2},d_{1}) \otimes \mathcal{Z}$.
\end{nprop}

\begin{nproof}
Set $B:= \mathrm{C}^{*}(d_{0},d_{1/2},d_{1})$, then $B \cong \mathcal{C}(Y)$ for some compact subset of the canonical 2-simplex $\Delta^{2} \subset \mathbb{R}^{3}$. With this identification, $C$ is isomorphic as a $\mathcal{C}(Y)$-algebra to a subalgebra  $\tilde{C}$ of $ \mathcal{C}(Y) \otimes \mathcal{Z} \otimes \mathcal{Z} \otimes \mathcal{Z}$. It is clear that for any $y \in Y$, the fibre $\tilde{C}_{y}$ is of the form $F_{1}\cdot F_{2} \cdot F_{3}$, where each $F_{i}$ is of the form $\mathcal{Z} \otimes \be_{\mathcal{Z}} \otimes \be_{\mathcal{Z}}$, $\be_{\mathcal{Z}} \otimes \mathcal{Z} \otimes \be_{\mathcal{Z}}$ or $\be_{\mathcal{Z}} \otimes \be_{\mathcal{Z}} \otimes \mathcal{Z}$. In any event, $A_{y}$ is isomorphic to a 1-, 2- or 3-fold tensor product of $\mathcal{Z}$ with itself, hence to  $\mathcal{Z}$. Since $\dim Y \le 2$, the main result of \cite{DadWin:trivial-fields} yields that $C \cong \mathcal{C}(Y) \otimes \mathcal{Z}$ as $\mathcal{C}(Y)$-algebras.
\end{nproof}
\en

\bn 
We now take up the task of proving Theorem~\ref{zstab}.

\label{crossed-products-Z-stable}
\begin{ntheorem}
Let $X$ be an infinite, compact, metrizable, finite-dimensional space with a minimal homeomorphism $\alpha$. Then, $\mathcal{C}(X) \rtimes_{\alpha} \mathbb{Z}$ is $\mathcal{Z}$-stable.
\end{ntheorem}

\begin{nproof}
Set
\[
A:= \mathcal{C}(X) \rtimes_{\alpha} \mathbb{Z}
\]
and let $u$ denote the unitary in $A$ implementing the action $\alpha$. Let $\epsilon >0$ and finite subsets $\mathcal{F} \subset A$ and $\mathcal{E} \subset \mathcal{Z}$ be given. In view of  Proposition~\ref{almost-central-D-stable} it will suffice to construct a unital $*$-homomorphism
\[
\zeta: \mathcal{Z} \to A_{\infty}
\]
such that
\begin{equation}
\label{w2}
\|[\zeta(z),a]\| < \epsilon
\end{equation}
for $z \in \mathcal{E}, a \in \mathcal{F}$. We may clearly assume that $\mathcal{E}$ and $\mathcal{F}$ consist of normalized elements, and that
\[
\mathcal{F} \subset \{\mathcal{C}(X) u^{l} \mid l=0, \ldots,k\}
\]
for some fixed $k \in \mathbb{N}$. 

Define elements of $\mathcal{C}([0,1])$ by
\[
c_{1}(t):= \left\{  
\begin{array}{ll}  1 & t=0\\
0 & t \ge \frac{1}{4}\\
\mbox{linear} & \mbox{else}
\end{array}  \right., \, 
c_{1/2}(t):= \left\{  
\begin{array}{ll}  0 & t=0,1\\
1 & \frac{1}{4} \le t \le \frac{3}{4}\\
\mbox{linear} & \mbox{else}
\end{array}  \right. , 
\]
\[
c_{0}(t):= \left\{  
\begin{array}{ll}  0 & t \le \frac{3}{4}\\
1 & t =1\\
\mbox{linear} & \mbox{else}
\end{array}  \right. , \, 
d_{1}(t):= \left\{  
\begin{array}{ll}  1 & t \le \frac{1}{4}\\
0 & t \ge \frac{1}{2}\\
\mbox{linear} & \mbox{else}
\end{array}  \right., 
\]
\[ 
d_{1/2}(t):= \left\{  
\begin{array}{ll}  0 & t \le \frac{1}{4}, t \ge \frac{3}{4}\\
1 & t =\frac{1}{2} \\
\mbox{linear} & \mbox{else}
\end{array}  \right. \mbox{ and }
d_{0}(t):= \left\{  
\begin{array}{ll}  0 & t \le \frac{1}{2}\\
1 & t \ge \frac{3}{4}\\
\mbox{linear} & \mbox{else}
\end{array}  \right. .
\]
The sets $\{c_1,c_{1/2},c_0\}$ and $\{d_1,d_{1/2},d_0\}$ thus form partitions of unity on $[0,1]$.  The rationale for our choice
of subscripts becomes clear just after (\ref{A-def}) below.

Set 
\[
\mathcal{G}:= \{ c_{i}^{\frac{1}{2k}}, d_{i} \mid i=0, 1/2, 1\}.
\]
Since $C$, defined as in Proposition~\ref{C-Z-trivial}, is trivial, there is a unital $*$-homomorphism 
\[
\tilde{\zeta}: \mathcal{Z} \to C.
\]  
We shall regard $C$ as a $\mathrm{C}^{*}$-subalgebra of 
\[
\tilde{C}:= \mathrm{C}^{*}(\mathcal{C}([0,1]) \otimes \be_{\mathcal{Z} \otimes \mathcal{Z} \otimes \mathcal{Z}}, C),
\]
but note that the elements of $\tilde{\zeta}(\mathcal{Z})$ are constant on $[0,\frac{1}{4}]$ and on $[\frac{3}{4},1]$.  

Next, use Proposition~\ref{C-gamma} with $\tilde{\zeta}(\mathcal{E})$ in place of $\mathcal{E}$, 
\[
\mathcal{S}:= (d_{0} \otimes \mathcal{Z} \otimes \be_{\mathcal{Z}} \otimes \be_{\mathcal{Z}}) \cup (d_{1/2} \otimes \be_{\mathcal{Z}} \otimes \mathcal{Z} \otimes \be_{\mathcal{Z}}) \cup (d_{1} \otimes \be_{\mathcal{Z}} \otimes \be_{\mathcal{Z}} \otimes \mathcal{Z}) 
\]
in place of $\mathcal{G}$ and with $\frac{\epsilon}{2}$ in place of $\theta$ to find $\beta >0$ and a finite
set $\mathcal{S}' \subset \mathcal{S}$.  The elements of $\mathcal{S}'$ are elementary tensors in which various elements of $\mathcal{Z}$ occur in the latter three tensor factors;  let $M$ denote the largest norm occurring among these elements (set $M=1$ if the said norm is zero).  Choose 
\[
0 < \eta < \frac{\epsilon}{6 k (k+1)}, \, \frac{\beta}{Mk}.
\]
Employ Proposition~\ref{function-commutator} to find $\delta>0$. Now for 
\[
x_{0} \neq x_{1} \in X,
\]
by Proposition~\ref{almost-central-h} there is 
\[
h \in \mathcal{C}(X) \subset A 
\]
positive and normalized such that 
\[
h(x_{0}) = 0, \, h(x_{1})= 1 \mbox{ and } \|[h,u]\| < \delta.
\]
For $i=0,1/2,1$ we set 
\[
\bar{c}_{i} := c_{i}(h) \mbox{ and } \bar{d}_{i} := d_{i}(h),
\]
and from our choice of $\delta$ and Proposition~\ref{function-commutator} we see that 
\[
\| [\bar{c}_{i}^{\frac{1}{2k}},u]\| , \, \| [\bar{d}_{i},u]\| < \eta
\]
for $i=0,1/2,1$. 

For convenience, we define
\begin{equation}\label{A-def}
A_{0}:= A_{\{x_{0}\}}, \, A_{1/2} := A_{\{x_{0},x_{1}\}} \mbox{ and } A_{1}:= A_{\{x_{1}\}}.
\end{equation}
With this notation we see that if $f \in \mathcal{C}([0,1])$ vanishes on $\{0\}$, then $u f(h) \in A_0$;  in particular, both 
$u \bar{c}_0$ and $u \bar{d}_0$ belong to $A_0$.  Similar containments hold for $A_1$ and $A_{1/2}$ when $f$ vanishes
on $\{1\}$ and $\{0,1\}$, respectively.  By Proposition~\ref{p-local-simplicity} the algebras of (\ref{A-def}) are simple with finite decomposition rank. By the main result of \cite{w2}, they are $\mathcal{Z}$-stable. But then, by \cite[Theorem~7.2.2]{Ror:encyc}, for $i=0,1/2,1$ there are $*$-homomorphisms 
\[
\varrho_{i}: A_{i} \otimes \mathcal{Z} \to (A_{i})_{\infty} \subset A_{\infty}
\]
such that
\[
\varrho_{i}|_{A_{i} \otimes \be_{\mathcal{Z}}} = \iota_{A_{i}}.
\]
Since $A_{1/2} \subset A_{i}$, we see that $\varrho_{1/2}(A_{1/2} \otimes \mathcal{Z}) \subset (A_{i})_{\infty}$ is a separable subset, whence by Proposition~\ref{diagonally-central} we may even assume that 
\begin{equation}
\label{w1}
[\varrho_{i}(\be_{A_{i}} \otimes \mathcal{Z}), \varrho_{1/2}(A_{1/2} \otimes \mathcal{Z})] = 0
\end{equation}
for $i=0,1$. We set
\[
\bar{\varrho}_{i}:= \rho_i \circ (\be_{A_{i}} \otimes \id_{\mathcal{Z}}): \mathcal{Z} \to (A_i)_\infty, \, i=0,1/2,1.
\]
Using \eqref{w1} and the fact that 
\[
\varrho_{0}(\bar{d}_{0} \otimes \mathcal{Z}) \perp \varrho_{1}(\bar{d}_{1} \otimes \mathcal{Z}),
\]
one checks that the assignment
\begin{eqnarray*}
f d_{0} \otimes z \otimes \be_{\mathcal{Z}} \otimes \be_{\mathcal{Z}}& \mapsto & \varrho_{0}(f(h) \bar{d}_{0} \otimes z) \\
f d_{1/2} \otimes \be_{\mathcal{Z}} \otimes z \otimes \be_{\mathcal{Z}}& \mapsto & \varrho_{1/2}(f(h) \bar{d}_{1/2} \otimes z) \\
f d_{1} \otimes \be_{\mathcal{Z}} \otimes \be_{\mathcal{Z}} \otimes z & \mapsto & \varrho_{1}(f(h) \bar{d}_{1} \otimes z) 
\end{eqnarray*}
for $f \in \mathcal{C}([0,1])$ and $z \in \mathcal{Z}$ extends to a unital $*$-homomorphism
\[
\gamma: \tilde{C} \to A_{\infty}.
\]
(For a sample calculation from this verification, let's see why
\[
[\gamma(f d_{0} \otimes z \otimes \be_{\mathcal{Z}} \otimes \be_{\mathcal{Z}}), \gamma(g d_{1/2} \otimes \be_{\mathcal{Z}} \otimes y \otimes \be_{\mathcal{Z}})]=0 
\]
for $f,g \in \mathcal{C}([0,1])$, and $z,y \in \mathcal{Z}$:
\begin{eqnarray*}
\lefteqn{ \gamma(f d_{0} \otimes z \otimes \be_{\mathcal{Z}} \otimes \be_{\mathcal{Z}}) \gamma(g d_{1/2} \otimes \be_{\mathcal{Z}} \otimes y \otimes \be_{\mathcal{Z}})} \\
& = & \varrho_0(f(h) \bar{d}_0 \otimes z) \varrho_{1/2}(g(h) \bar{d}_{1/2} \otimes y) \\
& = & \varrho_0(f(h) \bar{d}_0 \otimes \be_{\mathcal{Z}}) \varrho_0(\be_{A_0} \otimes z) \varrho_{1/2}(g(h) 
\bar{d}_{1/2} \otimes y)\\
& \stackrel{(\ref{w1})}{=} & (f(h) \bar{d}_0) \varrho_{1/2}(g(h) \bar{d}_{1/2} \otimes y) \varrho_0(\be_{A_0} \otimes z) \\
& = & (f(h) \bar{d}_0 g(h) \bar{d}_{1/2}) \varrho_{1/2}( \be_{A_{1/2}} \otimes y) \varrho_0(\be_{A_0} \otimes z) \\
& = & \varrho_{1/2}(\be_{A_{1/2}} \otimes y) (f(h) \bar{d}_0)( g(h) \bar{d}_{1/2})
\varrho_0(\be_{A_0} \otimes z)  \\
& = & \varrho_{1/2}(\be_{A_{1/2}} \otimes y) (g(h) \bar{d}_{1/2})(f(h) \bar{d}_0)
\varrho_0(\be_{A_0} \otimes z)  \\
& = & \varrho_{1/2}(g(h) \bar{d}_{1/2} \otimes y) \varrho_0(f(h) \bar{d}_0 \otimes z) \\
& = & \gamma(g d_{1/2} \otimes \be_{\mathcal{Z}} \otimes y \otimes \be_{\mathcal{Z}}) \gamma(f d_{0} \otimes z \otimes 
\be_{\mathcal{Z}} \otimes \be_{\mathcal{Z}}),
\end{eqnarray*}
as required.  Note that despite its length, our calculation uses only (\ref{w1}), the properties of $*$-homomorphisms on
tensor products, and the fact that the $\mathrm{C}^*$-algebra generated by $h$ commutes with everything of the form $\varrho_i(
\be_{A_i} \otimes z)$, where $i= 0,1/2,1$.)
We denote the restriction of $\gamma$ to $C$ also by $\gamma$, and so obtain a unital $*$-homomorphism
\[
\zeta:= \gamma \circ \tilde{\zeta} : \mathcal{Z} \to A_{\infty}.
\]

Our goal is to prove that for any $fu^l \in \mathcal{F}$ and $z \in \mathcal{E}$, we have
\[
\|[fu^l,\zeta(z)]\| < \epsilon.
\]
For $i=0,1/2,1$ we have 
\[
\| [ \bar{c}_{i}^{\frac{1}{2k}},u]\| < \eta,
\]
whence
\begin{eqnarray*}
\| \bar{c}_{i}f u^{l} - \bar{c}_{i}^{\frac{2k-l}{2k}} f (u \bar{c}_{i}^{\frac{1}{2k}})^{l} \| & < & l \eta + (l-1) \eta + \ldots + \eta \\
& \le & \frac{k(k+1)}{2} \eta \\
& \le & \frac{\epsilon}{12}.
\end{eqnarray*}
Set
\[
\tilde{c}_{i,l,f} := \bar{c}_{i}^{\frac{2k-l}{2k}} f (u \bar{c}_{i}^{\frac{1}{2k}})^{l}
\]
and note that 
\[
\tilde{c}_{i,l,f} \in \overline{\bar{c}_{i} A_{i} \bar{c}_{i}}.
\]
Now
\begin{eqnarray*}
\| [\zeta(z),fu^{l}] \| & = & \| [\zeta(z), (\bar{c}_{0} + \bar{c}_{1/2} + \bar{c}_{1}) f u^{l}] \| \\
& \le & \| [ \zeta(z), \tilde{c}_{0,l,f} + \tilde{c}_{1/2,l,f} + \tilde{c}_{1,l,f}] \| + 6 \frac{\epsilon}{12} \\
& \leq & \|[\zeta(z), \tilde{c}_{0,l,f}] \| + \|[\zeta(z), \tilde{c}_{1/2,l,f}] \| + \|[\zeta(z), \tilde{c}_{1,l,f}] \| + \frac{\epsilon}{2}.
\end{eqnarray*}
We will prove that the first and third terms of the last line are zero, while the second is strictly less than $\epsilon/2$. This will
complete the proof of the theorem.

\noindent
{\sc Case I.}  We prove that $[\zeta(z),\tilde{c}_{0,l,f}] = [\zeta(z),\tilde{c}_{1,l,f}] = 0$.   
For $i=0,1$ and $z \in \mathcal{Z}$, we have 
\begin{eqnarray*}
\bar{d}_{i} \zeta(z) & = & \bar{d}_{i} \gamma \tilde{\zeta}(z) \\
& = & \varrho_{i}(\bar{d}_{i} \otimes \be_{\mathcal{Z}}) \gamma (\tilde{\zeta}(z)) \\
& = & \gamma(d_{i} \otimes \be_{\mathcal{Z}} \otimes \be_{\mathcal{Z}} \otimes \be_{\mathcal{Z}}) \gamma (\tilde{\zeta}(z)) \\
& \subset & \gamma(d_{i} \cdot C).
\end{eqnarray*}
It follows that 
\[
\bar{d}_{0} \zeta(z) \bar{c}_{0}^{\frac{1}{2}} = \gamma(d_{0} c_{0}^{\frac{1}{2}} \otimes \tilde{\zeta}(z)(1) \otimes \be_{\mathcal{Z}} \otimes \be_{\mathcal{Z}})
\]
and
\[
\bar{d}_{1} \zeta(z) \bar{c}_{1}^{\frac{1}{2}} = \gamma(d_{1} c_{1}^{\frac{1}{2}} \otimes \be_{\mathcal{Z}} \otimes \be_{\mathcal{Z}}  \otimes \tilde{\zeta}(z)(0) ),
\]
whence, setting $\bar{i} = i+1 \ (\mathrm{mod} \ 2)$,
\begin{eqnarray*}
\bar{d}_{i} \zeta(z) \bar{c}_{i}^{\frac{1}{2}} & = & \varrho_{i}(\bar{d}_{i} \bar{c}_{i}^{\frac{1}{2}} \otimes \be_{\mathcal{Z}}) \varrho_{i}(\be_{A} \otimes \tilde{\zeta}(z)(\bar{i})) \\
& = &  \varrho_{i}( \bar{c}_{i}^{\frac{1}{2}} \otimes \be_{\mathcal{Z}}) \varrho_{i}(\be_{A} \otimes \tilde{\zeta}(z)(\bar{i})) \\
& = &  \bar{c}_{i}^{\frac{1}{2}} \varrho_{i}(\be_{A} \otimes \tilde{\zeta}(z)(\bar{i})) \\
& = &  \varrho_{i}(\be_{A} \otimes \tilde{\zeta}(z)(\bar{i})) \bar{c}_i^{\frac{1}{2}}\\
& = & \varrho_i(\bar{c}_i^{\frac{1}{2}} \otimes \tilde{\zeta}(z)(\bar{i})) \\
& = & \varrho_i(\bar{c}_i^{\frac{1}{2}} \otimes \tilde{\zeta}(z)) \\
& = &  \bar{c}_{i}^{\frac{1}{2}} \zeta(z) \\
& = & \bar{c}_i^{\frac{1}{2}} \bar{d}_i \zeta(z)
\end{eqnarray*}
for $i=0,1$. (To get from the fifth line to the sixth, recall that $\tilde{\zeta}(z)$ is constant on the support of $\bar{c}_i$.) Employing more than one of the equalities from the string immediately above, we see that for $i=0,1$ and $a \in A_{i}$
\begin{eqnarray*}
\bar{d}_{i} \zeta(z) \bar{c}_{i} a \bar{c}_{i} & = & \bar{d}_{i} \zeta(z) \bar{c}_{i}^{\frac{1}{2}} (\bar{c}_{i}^{\frac{1}{2}} a \bar{c}_{i})\\
& = &  \varrho_{i}(\be_{A} \otimes \tilde{\zeta}(z)(\bar{i})) \bar{c}_{i} a \bar{c}_{i} \\
& = &   \bar{c}_{i} a \bar{c}_{i} \varrho_{i}(\be_{A} \otimes \tilde{\zeta}(z)(\bar{i}))\\
& = &  (\bar{c}_{i} a \bar{c}_{i}^{\frac{1}{2}}) \bar{c}_{i}^{\frac{1}{2}} \bar{d}_{i} \zeta(z)\\
& = &  \bar{c}_{i} a \bar{c}_{i} \bar{d}_{i} \zeta(z),
\end{eqnarray*}
whence 
\[
[\bar{d}_{i} \zeta(\mathcal{Z}), \overline{\bar{c}_{i} A_{i} \bar{c}_{i}}] = 0
\]
for $i=0,1$. Since
\[
\bar{d}_{i} \zeta(z) = \zeta(z) \bar{d}_{i}  \mbox{ and } \bar{d}_{i}\bar{c}_{i} = \bar{c}_{i}
\]
for $i=0,1$ and $z \in \mathcal{Z}$, we have
\begin{eqnarray*}
[\zeta(\mathcal{Z}), \overline{\bar{c}_{i} A_{i} \bar{c}_{i}}] & = & [(\bar{d}_i + \bar{d}_{1/2} + \bar{d}_{\bar{i}}) \zeta(z), 
\overline{\bar{c}_{i} A_{i} \bar{c}_{i}}]  \\
& = & [\bar{d}_i \zeta(z), \overline{\bar{c}_{i} A_{i} \bar{c}_{i}}]  \\
& = & 0,
\end{eqnarray*}
as required.

\vspace{2mm}
\noindent
{\sc Case II.}  We prove that $\|[\zeta(z),\tilde{c}_{1/2,l,f}]\| < \epsilon/2$.
For any
\[
d_0 \otimes y \otimes 1_{\mathcal{Z}} \otimes 1_{\mathcal{Z}} \in \mathcal{S}'
\]
we have
\begin{eqnarray*}
\|[\tilde{c}_{1/2,l,f}, \gamma(d_0 \otimes y \otimes 1_{\mathcal{Z}} \otimes 1_{\mathcal{Z}})]\| & = & \| [\tilde{c}_{1/2,l,f} , \rho_0(\bar{d}_{i} \otimes y)] \| \\
& = & \| [ \tilde{c}_{1/2,l,f}, \rho_0(\bar{d}_0 \otimes 1_\mathcal{Z}) \rho_0(1_{A_i} \otimes y)] \| \\
& \leq & \| [\tilde{c}_{1/2,l,f}, \bar{d}_i]\| \cdot \| \rho_0(1_{A_i} \otimes y) \| \\
& \leq & \| [\tilde{c}_{1/2,l,f}, \bar{d}_i]\| \cdot M \\
& \le & M l \| [u, \bar{d}_{i}] \| \\
& < & M l \eta \\
& \le & \beta. 
\end{eqnarray*}
Similar calculations show that 
\[
\|[\tilde{c}_{1/2,l,f}, \gamma(s) ] \| < \beta, \ \forall s \in \mathcal{S}'.
\]
It then follows from our choice of $\beta$ and Proposition~\ref{C-gamma} that
\[
\| [\tilde{c}_{1/2,l,f}, \zeta(z)] \| < \frac{\epsilon}{2},
\]
as required.
\end{nproof}
\en

\section{Classification by $\mathrm{K}$-theory:  Elliott's conjecture}

\bn
Our proof of Theorem \ref{main} will be an application of a general classification result, due to the second named author (Theorem 7.1 of \cite{w3}).  To avoid the task of deriving a succinct and useful form of this theorem for our purpose, we state here a generalization of \cite[Theorem 7.1]{w3} due to Lin and Niu (\cite{ln}).
First we fix some notation:  if $l$ is a prime, then set
\[
\mathfrak{U}_l := \bigotimes_{i=1}^\infty \mathrm{M}_l.
\]

\begin{ntheorem}\cite{w3, ln}
Let $A,B$ be unital simple separable nuclear $\mathrm{C}^*$-algebras which absorb the Jiang-Su algebra $\mathcal{Z}$ tensorially.  Suppose
that for any prime $l$, the tensor products $ \mathfrak{U}_l \otimes A$ and $ \mathfrak{U}_l  \otimes B$ have tracial rank zero and satisfy the universal coefficient theorem.  Also suppose that there is a graded ordered isomorphism $\phi:\mathrm{K}_*(A) \to \mathrm{K}_*(B)$.  It follows that there is a $*$-isomorphism $\Phi:A \to B$ inducing $\phi$.
\end{ntheorem}
\en

\bn
\label{main-proof-reduction}
The crossed products considered in Theorem \ref{main} are simple, unital, separable, nuclear $\mathrm{C}^*$-algebras satisfying the universal coefficient theorem, and were shown in Theorem \ref{crossed-products-Z-stable} to be $\mathcal{Z}$-stable.  Thus, we see that the proof of Theorem \ref{main} is reduced to the problem of showing that for any prime $l$ and $\mathrm{C}^*$-algebra $A$ as in the said theorem, the tensor product $\mathfrak{U}_l \otimes A$ has tracial rank zero. 

This task will be accomplished in two steps: First, we will employ known results to show that, under the hypotheses of Theorem~\ref{main} and with notation as in \ref{A-Y-def}, algebras of the form $\mathfrak{U}_{l} \otimes A_{\{y\}}$ are TAF. This is the content of Proposition~\ref{B_y-tr0} below. 

Second, we need to conclude that if  $\mathfrak{U}_{l} \otimes A_{\{y\}}$ is TAF, then so is $\mathfrak{U}_{l} \otimes (\mathcal{C}(X) \rtimes_{\alpha} \mathbb{Z})$. This will be a special case of a result by  Strung and the second named author,  Theorem~\ref{TAS}; we will outline the proof at least in the special case we need, cf.\ Proposition~\ref{B-tr0} below. 
\en

%
%

\bn
So, let us turn to the first step of \ref{main-proof-reduction}.

\begin{nprop}\label{B_y-tr0}
Let $X$ be an infinite, compact, metrizable space of finite covering dimension, and $\alpha:X \to X$ a minimal homeomorphism.
Suppose further that the projections of $\mathcal{C}(X) \rtimes_\alpha \mathbb{Z}$ separate traces.  It follows that for any prime $l$ and any $y \in X$, the $\mathrm{C}^*$-algebra $B_{y}:= \mathfrak{U}_l \otimes A_{\{y\}}$
has tracial rank zero.
\end{nprop}

\begin{nproof}
By Proposition \ref{p-local-simplicity}, $A_{\{y\}}$ is a simple, $\mathcal{Z}$-stable $\mathrm{C}^*$-algebra with decomposition rank at most $d:=\mathrm{dim}(X)$.  The properties of the decomposition rank, established in \cite{kw}, include the following:  
decomposition rank is insensitive to taking tensor products with full matrix algebras over $\mathbb{C}$; decomposition rank
is lower semicontinuous with respect to inductive limits.  Since $\mathfrak{U}_l$ is an inductive limit of full matrix algebras over $\mathbb{C}$, we conclude that $B_y$ has decomposition rank at most $d$.  The inclusion $A_{\{y\}} \hookrightarrow A$
induces isomorphisms $\gamma_1:\mathrm{K}_0(A_{\{y\}}) \to \mathrm{K}_0(A)$ and $\gamma_2:\mathrm{T}(A) \to \mathrm{T}(A_{\{y\}})$
(see Theorem 4.1 (3) of \cite{p1} and \cite[Proposition 16]{ql}, respectively---a sketched proof of the latter result can be found in
\cite{qlp}).  It follows that projections separate traces in $A_{\{y\}}$.  Using the fact that $B_y$ is unital, simple, and absorbs $\mathfrak{U}_l$ tensorially, we may appeal to  \cite{r1} to conclude that $B_y$ has real rank zero.  We have now collected the hypotheses of \cite[Theorem 4.1]{w1}.
We conclude that $B_y$ has tracial rank zero, as desired.  
\end{nproof}
\en

This settles the first step of \ref{main-proof-reduction} above. To take care of the second step, we state the main technical result of \cite{StrWin:TAS}.  It refers to the concept of TA$\mathcal{S}$ algebras (where $\mathcal{S}$ is some class of separable, unital $\mathrm{C}^{*}$-algebras), which generalizes that of TAF algebras and was studied in detail in \cite{EllNiu:TAS}.

\bn
\label{TAS} 
\begin{ntheorem}
Let $X$ be an infinite, compact, metrizable space, and $\alpha: X \to X$ a minimal homeomorphism. Let $y \in X$, let $l$ be a prime, and suppose that $\mathfrak{U}_{l} \otimes A_{\{y\}}$ is TA$\mathcal{S}$, where $\mathcal{S}$ is some class  of unital semiprojective $\mathrm{C}^{*}$-algebras which is closed under taking quotients.  Then, $\mathfrak{U}_{l} \otimes (\mathcal{C}(X) \rtimes_{\alpha} \mathbb{Z})$ is TA$\mathcal{S}$ as well.
\end{ntheorem}

When we take $\mathcal{S}$ to be the class of finite dimensional $\mathrm{C}^{*}$-algebras, the theorem indeed settles the second step of \ref{main-proof-reduction} above. Since we will only need this special case, for the convenience of the reader we briefly outline its proof  in the sequel (cf.\ Proposition~\ref{B-tr0} below). We refer to \cite{StrWin:TAS} for full details, and for a discussion of the more general Theorem~\ref{TAS} and its implications for the classification program. 
\en

\bn
The key step in the proof of Theorem~\ref{TAS} is the following generalization of \cite[Lemma~4.2]{lp}, which will be derived  in \cite{StrWin:TAS}; here we only need it in the case of real rank zero. 

\begin{nlemma}\label{Bystructure}
Let $X$ be an infinite, compact, metrizable space, and $\alpha:X \to X$ a minimal homeomorphism.  Let $y \in X$, let $l$ be a prime and set $B_y = \mathfrak{U}_l 
\otimes A_{\{y\}}$.  It follows that for any $\epsilon>0$ and any finite subset 
$\mathcal{F}$ of $B:= \mathfrak{U}_l \otimes (\mathcal{C}(X) \rtimes_\alpha \mathbb{Z})$, there is a projection $p \in B_y$ such that the following 
statements hold:
\begin{enumerate}
\item[(i)] $\|pa - ap\| < \epsilon$ for all $a \in \mathcal{F}$;
\item[(ii)] $\mathrm{dist}(pap,pB_yp) < \epsilon$ for all $a \in \mathcal{F}$;
\item[(iii)] $\tau(1-p) < \epsilon$ for all $\tau \in \mathrm{T}(B)$. 
\end{enumerate}
\end{nlemma}

The proof is a modification of that of \cite[Lemma~4.2]{lp}. 
In fact, when we replace $\mathfrak{U}_{l}$ by $\mathbb{C}$ in the above (i.e., we take $l=1$), and assume in addition that $A_{\{y\}}$ has real rank zero and stable rank one, the lemma is exactly \cite[Lemma~4.2]{lp}. Now if  $\mathfrak{U}_{l}$ is nontrivial, and if $B_{y}$ has real rank zero, then the proof of  \cite[Lemma~4.2]{lp} carries over almost verbatim, at least in the case where $\mathcal{F}$ is a subset of $\be_{\mathfrak{U}_{l}} \otimes A_{\{y\}}$. The case of general $\mathcal{F}$ can be reduced to the situation where $\mathcal{F}$ is in $\be_{\mathfrak{U}_{l}} \otimes \mathfrak{U}_{l} \otimes \be_{A_{\{y\}}} \cup \be_{\mathfrak{U}_{l}} \otimes \be_{\mathfrak{U}_{l}} \otimes A_{\{y\}}$ (using the fact that $\mathfrak{U}_{l} \cong \mathfrak{U}_{l} \otimes \mathfrak{U}_{l}$); the proof in this situation is essentially the same as in the first case. 
\en

\bn
Next, we recall another result from \cite{lp}.


\begin{nlemma}\cite[Lemma 4.4]{lp}\label{lp4.4}
Let $C$ be a simple unital $\mathrm{C}^*$-algebra.  Suppose that for every finite subset $\mathcal{F} \subset C$, every $\epsilon > 0$, and
every nonzero positive element $b \in C$, there exists a projection $p \in C$ and a simple unital subalgebra $D \subset pCp$ with
tracial rank zero such that:
\begin{enumerate}
\item[(i)]  $\|pa-ap\| < \epsilon$ for all $a \in \mathcal{F}$;
\item[(ii)] $\mathrm{dist}(pap,D) < \epsilon$ for all $a \in \mathcal{F}$;
\item[(iii)] $p$ is Murray-von Neumann equivalent to a projection in $\overline{bCb}$.
\end{enumerate}
It follows that $A$ has tracial rank zero.
\end{nlemma}
\en

\bn
Using essentially the same technique as in \cite[Theorem~4.5]{lp}, we may now use the two preceding lemmas to prove Theorem~\ref{TAS},  at least in the case where $\mathcal{S}$ is the class of finite-dimensional $\mathrm{C}^{*}$-algebras. 

\begin{nprop}\label{B-tr0}
Let $X$ be an infinite, compact, metrizable space of finite covering dimension, and $\alpha:X \to X$ a minimal homeomorphism.
Suppose further that, for some prime $l$ and $y \in X$,  the $\mathrm{C}^*$-algebra $B_{y}:= \mathfrak{U}_l \otimes  A_{\{y\}}$ has tracial rank zero.  It follows that the $\mathrm{C}^*$-algebra $B:= \mathfrak{U}_l \otimes (\mathcal{C}(X) \rtimes_\alpha \mathbb{Z})$
has tracial rank zero.
\end{nprop}

\begin{nproof}
Since $B_y$ has tracial rank zero, so too does $pB_yp$ for any projection $p \in B_y$.  Using this fact and 
Lemma \ref{Bystructure} we may apply \cite[Lemma 4.4]{lp} (stated as Lemma \ref{lp4.4} above) with $C$ replaced by
$B:= \mathfrak{U}_l \otimes (\mathrm{C}(X) \rtimes_\alpha \mathbb{Z})$ and $D$ replaced by our $pB_yp$, where $p \in B_y$
is the projection provided by the conclusion of Lemma \ref{Bystructure}.  To conclude that $B$
has tracial rank zero, we need to show that the projection $p$ in the conclusion of Lemma \ref{Bystructure} is Murray-von Neumann
equivalent to a projection in $\overline{bBb}$
(we assume that $\epsilon$, $\mathcal{F}$, and $b$ as in the hypotheses of Lemma
\ref{lp4.4} are given).  By shrinking $\epsilon$, however, the existence of this projection follows from the second part of 
\cite[Proposition 3.8]{pt}---this proposition applies since $B$ is unital, simple, and $\mathcal{Z}$-stable.  We conclude, 
finally, that $B$ has tracial rank zero, as desired.
\end{nproof}
\en

We have thus completed the two steps of \ref{main-proof-reduction}, hence the proof of Theorem \ref{main}.

\end{document}